# Simpler near-optimal controllers through direct supervision

**Douglas B. Tweed**

*Departments of Physiology and Medicine, University of Toronto, Toronto, Ontario M5S 1A8, Canada; Centre for Vision Research, York University, Toronto, Ontario M3J 1P3, Canada*

**Abstract**

The method of generalized Hamilton-Jacobi-Bellman equations (GHJB) is a powerful way of creating near-optimal controllers by learning. It is based on the fact that if we have a feedback controller, and we learn to compute the gradient $\nabla J$ of its cost-to-go function, then we can use that gradient to define a better controller. We can then use the new controller's $\nabla J$ to define a still-better controller, and so on. Here I point out that GHJB works *indirectly* in the sense that it doesn't learn the best approximation to $\nabla J$ but instead learns the time derivative $\dot{J}$, and infers $\nabla J$ from that. I show that we can get simpler and lower-cost controllers by learning $\nabla J$ directly. To do this, we need teaching signals that report $\nabla J(\boldsymbol{x})$ for a varied set of states $\boldsymbol{x}$. I show how to obtain these signals, using the GHJB equation to calculate one component of $\nabla J(\boldsymbol{x})$ — the one parallel with $\dot{\boldsymbol{x}}$ — and computing all the other components by backward-in-time integration, using a formula similar to the Euler-Lagrange equation. I then compare this direct algorithm with GHJB on 2 test problems.

## 1 Introduction

Optimal control is computationally so demanding that it is probably unattainable except for low-dimensional tasks (Bellman 1961). The next best thing may be *near*-optimal control, where we compute a sequence of better and better controllers, moving ever closer to the optimal one (Sutton and Barto 1998; Sutton, McAllester et al. 2000). There are many algorithms for this purpose, but here I am concerned with methods involving continuous state spaces and time, and of these the most efficient, in terms of time and memory, is probably the method of *generalized Hamilton-Jacobi-Bellman* (GHJB) equations (Saridis and Lee 1979; Beard, Saridis et al. 1997; Lyshevski 1998; Abu-Khalaf and Lewis 2004; Abu-Khalaf and Lewis 2005). I will show that this method is in an important sense *indirect*, and that it can be modified to give better controllers with fewer adjustable parameters by using a more direct form of supervised learning.

## 2 Near-optimal control

### 2.1 Cost-to-go

Most methods of near-optimal control depend on the idea of *cost-to-go* functions, or some variant such as action-value or Bellman functions (Sutton and Barto 1998; Todorov 2004; Abu-Khalaf and Lewis 2005). To define the cost-to-go, we start with the *plant* equation

$$\dot{x} = \overline{f}(x, u(x)) = f(x) + G(x)u(x) \qquad (2.1)$$

which I will usually write as $\dot{x} = f + Gu$. Here $x$ is the plant state and $u$ is the control signal or *command*, which is computed from $x$ by a function called a *control law* or *controller*. In (2.1) the plant is affine in $u$, but the method can be generalized to the non-affine case by the usual method of treating $\dot{u}$ as the command and regarding $u$ as part of an augmented state.

The aim of optimal control is to find the controller that minimizes some *cost*

$$C = \int_0^\infty L(x, u) dt \qquad (2.2)$$

where $L$ is a nonnegative scalar called the *cost rate* or *loss*, e.g. $L$ might be the distance from a person's or a robot's gaze point or hand to a target. (Throughout this paper I will assume an infinite time horizon, as in (2.2), though the ideas can be modified to deal with finite time spans. Further, (2.2) has the loss depend only on the state and command, but it can be generalized; e.g. to have the loss depend on a target state $x^*$ we can define a new, augmented state vector that incorporates $x^*$.)

In near-optimal control we construct a sequence of better and better controllers $u^{(1)}$, $u^{(2)}$, $u^{(3)}$ and so on. Given any one controller $u^{(n)}$ we define its *cost-to-go* $J^{(n)}$ like this:

$$J^{(n)}(x) = \int_t^\infty L(x, u^{(n)}(x)) dt \tag{2.3}$$

i.e., $J^{(n)}(x)$ is the total cost we will accrue from now on if we are currently (at time $t$) in state $x$ and we choose our commands using controller $u^{(n)}$.

## 2.2  Building better controllers with the gradient of the cost-to-go

If we know $J^{(n)}$ for one controller $u^{(n)}$ then we can define a better controller $u^{(n+1)}$ [ref]. Specifically, we define $u^{(n+1)}$ to be the minimizer of $\dot{J}^{(n)} + L$, i.e.

$$\begin{aligned}
\forall x \; u^{(n+1)}(x) &= \arg\min_{u \in U} \dot{J}^{(n)} + L \\
&= \arg\min_{u \in U} \nabla J^{(n)} \dot{x} + L \\
&= \arg\min_{u \in U} \nabla J^{(n)}(x)[f(x) + G(x)u] + L(x, u)
\end{aligned} \tag{2.4}$$

where $U$ is the set of permissible commands and $\nabla J^{(n)}$ is the gradient of $J^{(n)}$ with respect to $x$, $dJ^{(n)}/dx$ (Figure 1). For instance if $L = x^T Q x + u^T u / 2$ and $u$ is unconstrained then the minimizer is

$$u^{(n+1)}(x) = -\left[\nabla J^{(n)}(x) G(x)\right]^T \tag{2.5}$$

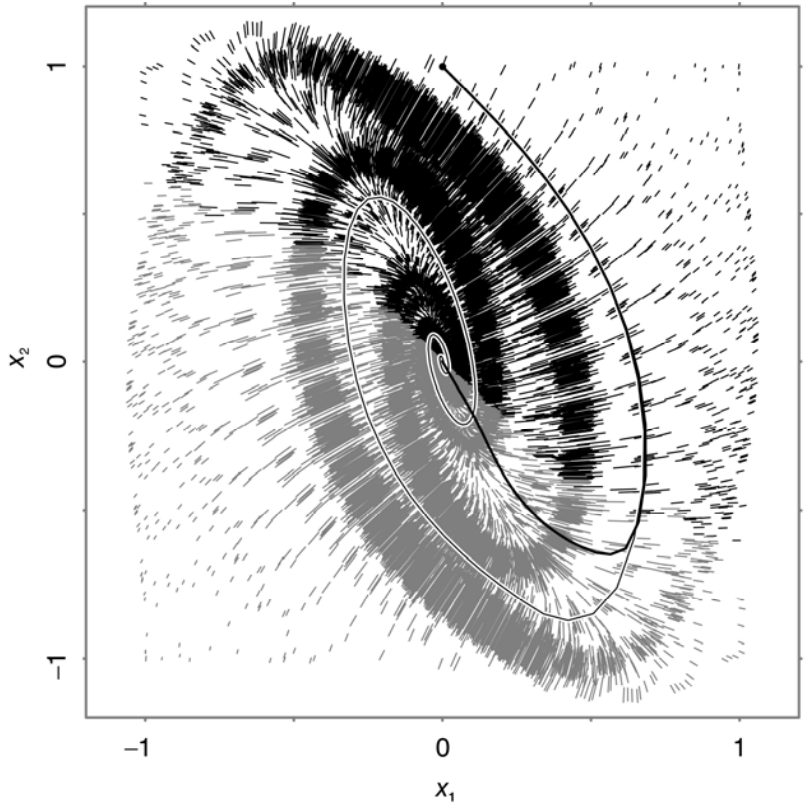

**Figure 1.** An example of the gradient field $\nabla J^{(1)}$ of an initial controller. Grey vectors are ones for which the product $\nabla J^{(1)}(x)G(x)$ is negative; for the black vectors the product is positive. The thin curve is the trajectory generated by $u^{(1)}$ on a test movement; the thick line is the trajectory of the improved controller $u^{(2)}$.

As (2.5) shows, we need the gradient, $\nabla J^{(n)}$, to create the new, improved controller $u^{(n+1)}$. And while this equation applies to just one specific loss function, it is clear from equation (2.4) that, whatever the loss, the formula for the minimizer will include $\nabla J^{(n)}$. So it is this gradient, rather than the cost-to-go itself, that is crucial.

We repeat the process, using the gradient of the new cost-to-go, $\nabla J^{(n+1)}$, to define a still-better controller $u^{(n+2)}$, and so on. If we can get an exact estimate of $\nabla J^{(n)}$ for each $n$, then our sequence is guaranteed to converge to the optimal controller $u^*$ (Abu-Khalaf and Lewis 2005). But in practice the estimates are never exact, and the main challenge is getting useful approximations to $\nabla J^{(n)}$.

## 2.3  Indirect approximation in GHJB

In GHJB, we learn to approximate $\nabla J^{(n)}$, but we do so indirectly. We differentiate (2.3) to get

$$\dot{J}^{(n)} = -L \qquad (2.6)$$

This equation tells us we can use $L$ as a teaching signal to train an approximator of $\dot{J}^{(n)}$ by supervised learning. The fastest and most reliable method, and the one adopted by Abu-Khalaf and Lewis (2005), is to use linear-in-the-parameters learning. First we define a *feature vector* $\theta(x)$ — a vector of nonlinear, continuously differentiable functions of $x$. Then we differentiate $\theta$ with respect to time to get $\dot{\theta}$, or equivalently we compute $\dot{\theta} = \nabla \theta \dot{x} = (d\theta/dx)\dot{x}$. And finally, we calculate the row vector of weights $w$ that provides the best approximation

$$\hat{\dot{J}} = w\dot{\theta} \qquad (2.7)$$

That is, we find the $w$ that minimizes the sum of squared approximation errors over some data set of $n$ examples:

$$SSE = \sum_{i=1}^{N} e^2 = \sum_{i=1}^{N} \left[ w\dot{\theta}(x_i) - \dot{J}^{(n)}(x_i) \right]^2 = \sum_{i=1}^{N} \left[ w\dot{\theta}(x_i) + L(x_i) \right]^2 \qquad (2.8)$$

The resulting optimal $w$ — the weight vector that yields the best approximation to $\dot{J}^{(n)}$ — I will call $w^{(n)}$.

With $w^{(n)}$ in hand, we have our approximation to $\dot{J}^{(n)}$, namely $w^{(n)}\dot{\theta}$, and also an approximation to the gradient $\nabla J^{(n)}$,

$$\nabla \hat{J}^{(n)} = w^{(n)} \nabla \theta = w^{(n)} \frac{d\theta}{dx} \qquad (2.9)$$

So this is the GHJB method: learn $\dot{J}^{(n)}$, derive an estimate of $\nabla J^{(n)}$, and use that to create $u^{(n+1)}$.

## 3  Equivalence-breaking

This technique of indirect learning succeeds if the approximation $\hat{\dot{J}}^{(n)}$ is very accurate. For instance in the extreme case, if our estimate of $\dot{J}^{(n)}$ is exactly right — that is, if $w^{(n)}\dot{\theta} = \dot{J}^{(n)}$ — then of course we have $\nabla J^{(n)} = w^{(n)} \nabla \theta$: our estimate of $\nabla J^{(n)}$ is also exact. So in this case, learning $\dot{J}^{(n)}$ is equivalent to learning $\nabla J^{(n)}$. But in practice our approximations are always imperfect, and the equivalence is broken: the vector $w^{(n)}$, obtained by minimizing (2.8), gives us the best approximation *to the time-derivative* $\dot{J}^{(n)}$ that is possible given the feature vector $\theta$ — the best approximation of the form $w\dot{\theta}$ — but this $w^{(n)}$ may yield an estimate of $\nabla J^{(n)}$ that is very far from optimal.

Here I suggest that we should instead minimize

$$SSE' = \sum_{i=1}^{N} e_i'^2 = \sum_{i=1}^{N} \left[ w \nabla \theta(x_i) - \nabla J^{(n)}(x_i) \right]^2 \qquad (3.1)$$

Because with this approach we arrive at a different weight vector, $w^{(n)'}$, which gives us the best approximation *to the gradient* $\nabla J^{(n)}$ of the form $w\nabla\theta$.

So the message is that, for any given set of features, we get better estimates of $\nabla J^{(n)}$, and therefore better controllers, if we learn $\nabla J^{(n)}$ directly, rather than inferring it from $\hat{\dot{J}}^{(n)}$, because then all our approximating power is focused on the quantity that matters — $\nabla J^{(n)}$, not $\dot{J}^{(n)}$. And for the same reason we can get equivalent or better controllers with lower-dimensional $\theta$, so the controllers are simpler in the sense that they involve fewer functions and adjustable parameters.

## 4 Learning the gradient

### 4.1 Backward integration

How can we learn $\nabla J^{(n)}$ directly? For any kind of supervised learning we need teaching signals, which are examples of the outputs of the function we are trying to approximate. In the GHJB method, for instance, the teaching signals are the scalars $\dot{J}^{(n)}(x_i)$, or in other words $-L(x_i, u^{(n)}(x_i))$, for many different state vectors $x_i$, as shown in (2.8).

To learn $\nabla J^{(n)}$ directly, then, we have to compute examples, $\nabla J^{(n)}(x_i)$. This we can do by integrating a differential equation, in much the same way that one computes costates $\lambda$ by integrating the Euler-Lagrange equation (Bryson and Ho 1975) backward in time:

$$\dot{\lambda} = -\frac{\partial L}{\partial x} - \lambda \frac{\partial \bar{f}}{\partial x} \qquad (4.1)$$

We can't use the Euler-Lagrange equation itself because its costates are not the teaching signals we need, though they are related: costates are derivatives of $J$ with respect to $x$, *holding constant the control trajectory $u(t)$*, whereas what we need are derivatives of $J$ with respect to $x$, *holding constant the feedback control law, $u^{(n)}$*. Obviously these are two different things: if we hold the feedback controller constant, then a small change in $x$ will not usually leave the control trajectory unchanged.

The equation we need is

$$\frac{d}{dt} \nabla J^{(n)} = -\frac{dL^{(n)}}{dx} - \nabla J^{(n)} \frac{d\bar{f}^{(n)}}{dx} = -\nabla L^{(n)} - \nabla J^{(n)} D\bar{f}^{(n)} \qquad (4.2)$$

which resembles Euler-Lagrange except that the partial derivatives of $L$ and $\bar{f}$ in (4.1) have become total derivatives of $L^{(n)}$ and $\bar{f}^{(n)}$, where

$$L^{(n)}(x) = L\left(x, u^{(n)}(x)\right)$$
$$\overline{f}^{(n)}(x) = \overline{f}\left(x, u^{(n)}(x)\right)$$
(4.3)

According to (2.3), we have to integrate (4.2) starting from $t = \infty$, but in practice we can shave off most of this time because all our controllers, $u^{(1)}$, $u^{(2)}$ and so on, are *stabilizing*, by which I mean they yield a finite total cost (2.2), right out to time infinity, for any movement in some region of state space we care about. We can choose $u^{(1)}$ with this property, and then the algorithm ensures that all subsequent controllers have it as well. Because the cost integral (2.2) approaches a finite limit, we know that beyond some finite time the controller will accumulate negligible further cost. For instance, if the movement is directed at some target state $x^*$ then we can run the controller until it brings $x$ very close to $x^*$ and then integrate (4.2) back from there to compute $\nabla J^{(n)}$ for use in the direct supervised learning rule (3.1).

**4.2 Variants**

I have said we should focus our approximating power on $\nabla J^{(n)}$ rather than on $\dot{J}^{(n)}$ because $\nabla J^{(n)}$ is what we need to construct the next controller, $u^{(n+1)}$. But that argument can be taken a step further: what we actually need is the product $\nabla J^{(n)} G$, as in (2.5), so we should really take this as our supervisor signal, and minimize

$$SSE'' = \sum_{i=1}^{N} e_i''^2 = \sum_{i=1}^{N} \left[ w \nabla \theta(x_i) G(x_i) - \nabla J^{(n)}(x_i) G(x_i) \right]^2 \quad (4.4)$$

In fact there are advantages to both approaches, (3.1) and (4.4) — briefly, (4.4) yields a slightly better approximation to $\nabla J^{(n)} G$, while (3.1) yields a much better estimate of $J^{(n)}$, which is useful for geometric integration of the plant equation (2.1). In this paper I will present simulations based on (4.4).

Because numerical integration is inexact, (4.2) will yield inexact values of $\nabla J^{(n)}(x)$, so it is often useful to exploit an additional source of information about that gradient, namely the GHJB equation (2.6), which specifies one component of $\nabla J^{(n)}(x)$ — the one along the state-velocity vector $\dot{x}$:

$$-L(x, u^{(n)}(x)) = \dot{J}^{(n)}(x) = \nabla J^{(n)}(x)\dot{x} \qquad (4.5)$$

In my examples in Section 5, I calculate estimates of $\nabla J^{(n)}(x)$ using (4.2) and then apply the smallest adjustment that makes them fit (4.5).

One final technicality: backward-in-time integration of (4.2) can be very sensitive to boundary conditions near the target state $x^*$. This sensitivity, together with the inaccuracies of numerical integration, means that if we perform a movement, running the plant equation forward in time until $x$ equals some state $x'$ near $x^*$, and then we integrate (4.2) backwards starting from $x'$, we won't retrace the state trajectory of the original, forward-in-time movement. As a result, we won't be evaluating our $L^{(n)}$ and $\bar{f}^{(n)}$, from (4.3), on the same $x$'s as in the original state trajectory. In fact, the backward-in-time trajectories are often confined to a small submanifold of state space, so we don't get a good sample of $x_i$'s for our learning rule (4.4). But we can solve this problem using "breadcrumbs": during the forward-in-time integration, we store a series of $x_i$'s along the state trajectory; then during the backward-in-time integration we use these stored values to guide us back along the same path, as illustrated in the m-file Jx_ODE in the Appendix.

## 5 Numerical examples

I will compare GHJB and direct supervision using two numerical examples presented by Abu-Khalaf and Lewis.

## 5.1 Nonlinear oscillator with constrained input

In Example 5.2 from their 2005 paper, Abu-Khalaf and Lewis consider a control task where the plant equation is

$$\dot{x}_1 = x_1 + x_2 - x_1(x_1^2 + x_2^2), \quad \dot{x}_2 = -x_1 + x_2 - x_2(x_1^2 + x_2^2) + u \quad (5.1)$$

and the loss is

$$L = \tanh(x^\top x) + 2u \tanh^{-1}(u) + \log(1 - u^2) \quad (5.2)$$

For their initial controller, Abu-Khalaf and Lewis choose a linear function saturated at 1 and $-1$:

$$u^{(1)}(x) = \operatorname{sat}_{-1}^{+1}(-5x_1 - 3x_2) \quad (5.3)$$

Their feature vector $\theta$ has 24 components, comprising all of the monomials in $x_1$ and $x_2$ of orders 2, 4, 6 or 8 — for instance $\theta_1(x) = x_1^2$, $\theta_2(x) = x_1 x_2$, $\theta_{22}(x) = x_1^2 x_2^6$. They train their controllers over the state-space region $-1 \leq x_1 \leq 1$, $-1 \leq x_2 \leq 1$, and they illustrate the results by plotting state and control trajectories versus time for a test movement from an initial state $x = (0, 1)^\top$ to $\mathbf{0}$, as performed by their $u^{(1)}$ and their final controller, obtained after the algorithm has converged.

Abu-Khalaf and Lewis say they need eighth-order monomial features to achieve convergence — sixth-order isn't enough. I confirm this finding in Figure 2, which plots the total costs, for the test movement, for different sequences of controllers. The thick grey line, representing a run of GHJB with 2-, 4-, 6- and 8th-order monomials, converges to a cost of roughly 4.6. But the *thin* grey line, a run of GHJB with just 2-, 4- and 6th-order monomials, or 15 features in all, veers far from optimality: the test costs for $u^{(4)}$ and all subsequent controllers are well above the range for my plot.

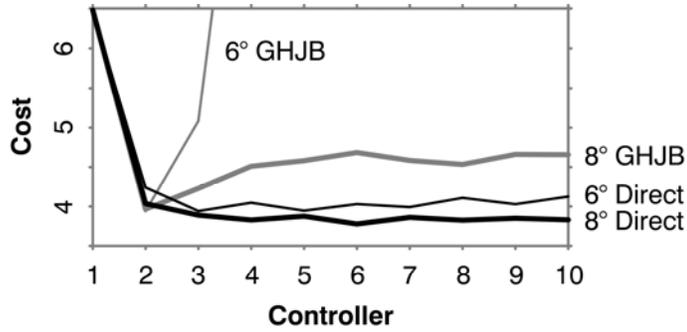

**Figure 2.** Tested on the control task from Example 5.2 of Abu-Khalaf and Lewis (2005), direct supervision (black) converges with 6th-order monomial features whereas GHJB (grey) needs 8th-order. But if we compare based on the best controller in each run, GHJB is about equally good with 6th or 8th-order features, and only 8th-order direct does better.

Direct supervision, on the other hand, converges with either feature vector, the 24- or the 15-component, as shown by the black lines in Figure 2. In other words, monomials to order 6 do suffice if all their approximating power is concentrated on the crucial function, $\nabla J^{(n)}$, as in the direct method, rather than on $J^{(n)}$, as in GHJB.

Convergence, though, isn't the central issue here, because it is neither necessary nor sufficient for good control. For instance, the "bad" sequence in Figure 2 (the thin grey line) does lead to some disastrous controllers, but only after producing one very good one: its $u^{(2)}$ achieves a cost of just 3.95 on the test movement (and there is nothing special about that movement — $u^{(2)}$ is about equally close to optimal for all other movements in the training region). Moreover convergence, when it does happen, need not bring us to the optimum. In Figure 2, the *thick* grey line descends to a cost of just 3.95 in round 2 but then deteriorates and converges at a much higher level (and again, there is nothing atypical about the test movement in Figure 2 — $u^{(2)}$ is better than its successors all-round). So if we judge based on the quality of the best controller in the sequence, rather than on convergence, then GHJB is about equally good with or without the 8th-order monomials, and 6th-order direct supervision is about the same. Only 8th-order direct supervision does better.

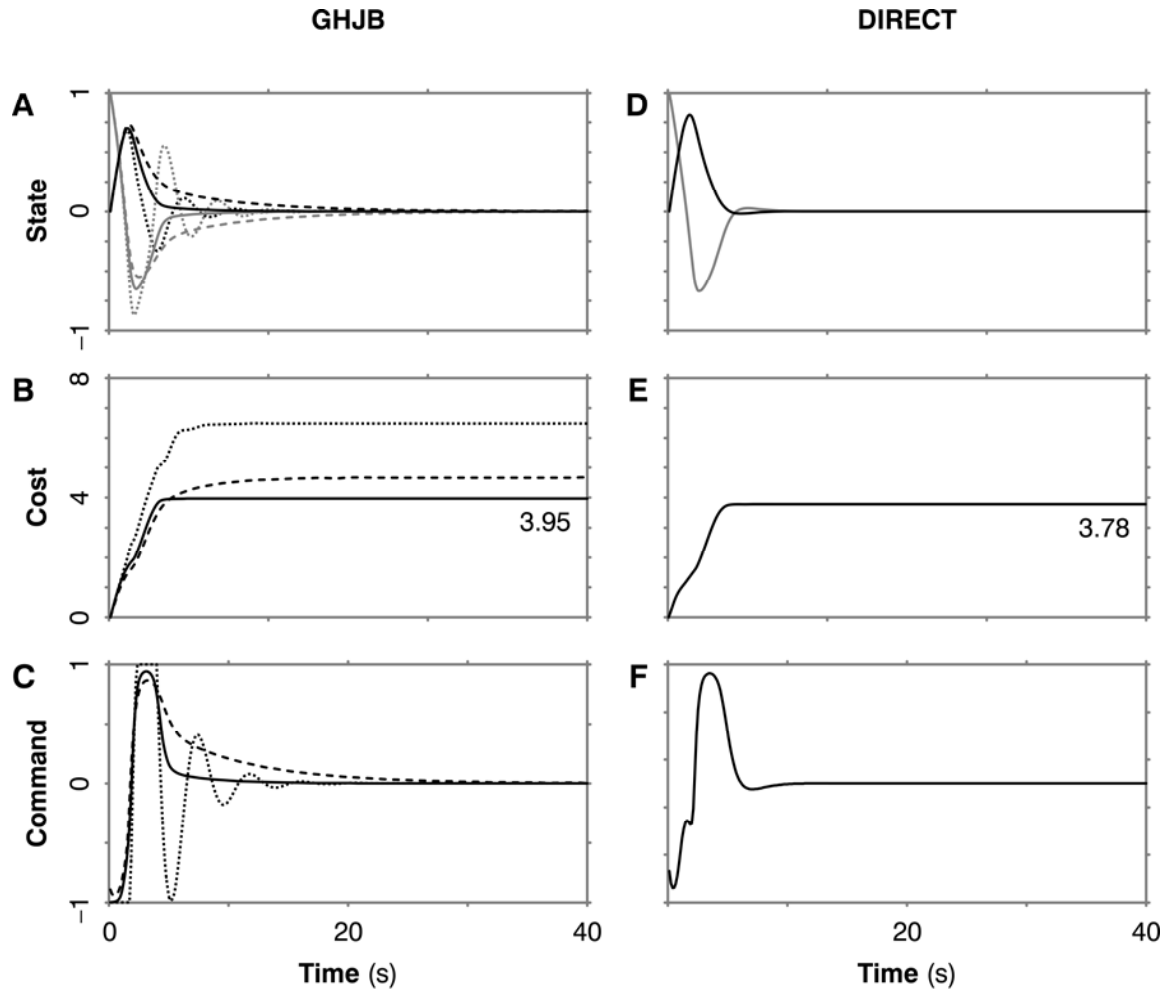

Figure 3. Comparing GHJB with direct supervision on the test movement from Example 5.2 of Abu-Khalaf and Lewis (2005).

Figure 3 shows trajectories for some of these controllers. The first column shows my recreation of Abu-Khalaf and Lewis's GHJB approach with all 24 features, up to eighth order. In Figure 3A, the dotted lines are the state trajectory ($x_1$ in black and $x_2$ in grey) for the test movement performed by the initial controller (5.3) — it matches the correspond-ing plot, Figure 5, in Abu-Khalaf and Lewis (2005). Dashed lines show the trajectories generated by the final controller, the limit of the sequence computed by GHJB; these curves match Abu-Khalaf and Lewis's "state trajectory for the nearly optimal control law" in their Figure 6. Finally, the solid lines in my Figure 3A shows the state trajectory for the *best* controller in the GHJB run, $u^{(2)}$. It accumulates a total cost, for the test movement, of just 3.95, whereas the final controller runs to about 4.6, as shown in 3B. Panel 3C shows the control trajectories for the same three controllers. Again, the dotted and dashed lines match the corresponding plots in Abu-Khalaf and Lewis.

In the second column of Figure 3 I show the trajectories for the direct supervision algorithm, also with features up to 8th order. The initial controller is the same as for GHJB, and the final controller doesn't differ much from the best one, so I plot only the best. Its cost is 3.78, about 4% lower than with GHJB.

**5.2 Random features**. My second comparison is similar to Example 5.3 from Abu-Khalaf and Lewis (2004), though not identical because some of the details are missing from that paper. The plant equation is

$$\dot{x}_1 = x_2, \ \dot{x}_2 = -x_2 + u \tag{5.4}$$

and the loss is

$$L(\boldsymbol{x}, u) = \tanh(\boldsymbol{x}^\mathsf{T} Q \boldsymbol{x}) + R\left[2u \tanh^{-1}(u) + \log(1 - u^2)\right] \tag{5.5}$$

where $Q = 100$ and $R = 1$. Abu-Khalaf and Lewis don't describe their $\boldsymbol{u}^{(1)}$, so I used a linear stabilizing controller,

$$u^{(1)}(\boldsymbol{x}) = -x_1 - x_2 \tag{5.6}$$

Regarding features, Abu-Khalaf and Lewis say only that they used "a smooth function of the order 35". I chose log-cosh functions

$$\boldsymbol{\theta} = \log(\cosh(\boldsymbol{W}_\theta \boldsymbol{x})) \tag{5.7}$$

where $W_\theta$ is a fixed matrix. And like Abu-Khalaf and Lewis, I trained the controllers over the state-space region $-0.5 \leq x_1 \leq 0.5$, $-0.5 \leq x_2 \leq 0.5$. The test movement I use for Figure 4 is the same as Abu-Khalaf and Lewis's, going from $x = (0.4, 0.4)^T$ to $\mathbf{0}$.

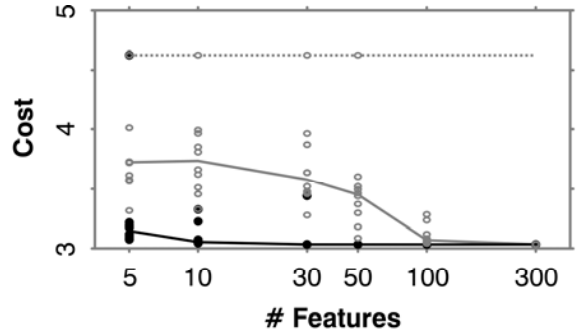

**Figure 4.** Comparing GHJB (grey) with direct supervision (black) on a more complex control task. Direct supervision achieves comparable performance to GHJB with about 10 times fewer features.

Figure 4 compares GHJB (the grey data) with direct supervision (black) for many runs with different numbers of features. On the abscissa is the number of features, and on the ordinate the cost of the test movement for the best controller in a any given run of the algorithm. For instance, there are 10 grey ellipses arranged in a vertical line above the tick mark for 50 on the abscissa. These represent the best costs achieved in 10 runs of GHJB with 50 log-cosh features, each run consisting of 19 rounds, or iterations, and therefore 20 controllers, $u^{(1)}$, $u^{(2)}$ … $u^{(20)}$. The matrix $W_\theta$, and therefore the feature functions, stay the same through all 19 rounds of any single run, but they differ from one run to the next because a new $W_\theta$ is chosen randomly at the start of each. Some of these random $W_\theta$ yield better features than others, and this is why the performance of the best controller varies from run to run. Near the top of the graph, the horizontal dotted grey line indicates the cost of the test movement for the *initial* controller, so ellipses on this line represent runs where none of the 20 controllers did better than $u^{(1)}$.

Also in a vertical line above this same tick mark (the one on the abscissa representing 50 features) are 10 black dots, clumped together near the bottom of the graph, representing the best costs achieved by direct supervision. These black dots lie below the grey ellipses, which means that direct supervision consistently finds better controllers than GHJB does when both algorithms work with 50 log-cosh features.

The solid lines mark median performance with different numbers of features — the grey line for GHJB, the black line for direct supervision. These lines show that direct supervision needs about 10 times fewer features than GHJB to achieve comparable performance. For instance direct supervision achieves a median cost within 1% of optimal with just 30 features, whereas GHJB needs 300 features to do as well. Even with just 5 features, direct supervision more often than not achieves costs within 4% of optimal. So as a rough estimate, if we run direct supervision on this problem 10 times with random sets of 5 features then we will find at least one controller that is at least this good with a probability of about $1 - 2^{-10}$, or 0.999.

I have also compared GHJB with direct supervision on Example 5.2 of Abu-Khalaf and Lewis (2004), but in that example both methods worked consistently with as few as 5 random log-cosh features. Direct supervision achieved costs that were 3% lower, but needed just as many features as GHJB. Presumably when GHJB gets by with so few features it leaves little room for improvement. But in more complex tasks such as the one in Figure 4, direct supervision can bring substantial savings.

The main drawback of direct supervision is that it learns more slowly than GHJB, at least in the cases I have tested. GHJB can collect training examples very quickly by sampling randomly from state space (Abu-Khalaf and Lewis 2005). Direct supervision, on the other hand, has to integrate the plant equation forward in time to a state $x'$ near the target $x^*$, then integrate (4.2) backward, as I have described. Of course each of these to-and-fro sweeps picks up a large number of training examples — one for each state $x_i$ along the path — but many of these states are close together, and they all lie on a single curve, so we need several such sweeps to get a representative sample of state space. In many tasks, though, getting a simpler controller may be worth the wait. And on complex problems, direct supervision may learn faster than GHJB, because learning speed depends on the number of adjustable parameters, and direct supervision reduces the required number.

Finally, direct supervision has other shortcomings which it shares with GHJB and most other algorithms for near-optimal control. It of course fails if we give it poor features or too few. Moreover, a set of features may suffice to fit $\nabla J^{(1)}$, in the first round of the algorithm, but then fail to fit $\nabla J^{(2)}$ because it has a more complicated shape. And sometimes we may get a good fit to $\nabla J^{(n)}$ but nevertheless a poor controller $u^{(n+1)}$ owing to domain effects — we learn $\nabla J^{(n)}$ over some region of state space but then $u^{(n+1)}$ carries us outside that region. Clearly we can solve all these problems by using a lot of features and a large domain, but more sophisticated variants might monitor the progress of the algorithm and adjust the features and training region as needed.

## 6 Conclusion

I have described an algorithm for near-optimal control that resembles the method of generalized Hamilton-Jacobi-Bellman equations but can yield equivalent or better controllers with fewer features. It does so by learning the least-mean-squares approximation to the gradient of the cost-to-go, $\nabla J^{(n)}$, rather than learning an approximation to $J^{(n)}$ and deriving an estimate of $\nabla J^{(n)}$ from that.

In assessing GHJB and direct supervision, I have argued that convergence should not be the main criterion, because it is neither necessary nor sufficient for finding good controllers. Instead I have focused on 2 other measures of success: the *best* controller obtained in any given run of the algorithm (whether that run later converges or diverges) and the me-dian performance of these best controllers when the algorithm is put through multiple runs with randomly chosen features.

**Acknowledgments**

I thank K. Fortney for comments on the manuscript. This work was funded by CIHR.

**Appendix**

The following is MATLAB code for direct supervision solving the control problem of Section 5.1 with 6th-order monomial features. The main m-file is the first, AK_2005_Ex_5_2_direct.m; the other 7 are called by it.

```
% AK_2005_Ex_5_2_direct.m
% Based on Abu-Khalaf & Lewis (2005), Numerical example 5.2.
% Direct supervision with 6th-order features.

clear all; clc

global G w_prev

% Initialize
nx = 2; nf = 15; nm = 100; nr = 5;   % # state dim'ns, features, rounds
ni = 10; nj = ni; np = ni*nj; nm_test = 1;   % # grid points, test
states
UCF = 0.1; wait = 40/UCF; % unit conversion factor for integration
G = UCF*[0; 1];
w = zeros(1, nf); w_prev = w;
b = 1;
X_TEST = b*(2*rand(nx, nm_test) - 1);   % col's are test states
X_TEST(:, nm_test) = [0; 1];

% Functions
u = @u1_AK52;   % initial controller
Du_Dx = @Du1_Dx_AK52;
f = @(x, ux) UCF*[x(1) + x(2) - x(1)*(x(1)^2 + x(2)^2); ...
                  -x(1) + x(2) - x(2)*(x(1)^2 + x(2)^2)] + G*ux;
df_dx = @(x) UCF*[1 - 3*x(1)^2 - x(2)^2, 1 - x(1)*2*x(2); ...
                  -1 - x(2)*2*x(1), 1 - x(1)^2 - 3*x(2)^2];
df_du = G;
F = @(x) f(x, u(x));
DF_Dx = @(x) df_dx(x) + df_du*Du_Dx(x);
L = @(x, ux) tanh(x'*x) + 2*ux*atanh(ux) + log(1 - (ux*ux));
L_ = @(x) L(x, u(x));
dL_dx = @(x) (1 - tanh(x'*x)^2)*2*x';
```

```
  dL_du = @(ux) 2*atanh(ux);
  DL_Dx = @(x) dL_dx(x) + dL_du(u(x))*Du_Dx(x);
  z = @z_AK52_246;   % features
  dz_dx = @dz_dx_AK52_246;

  % Compute costs of initial controller
  C_sum = 0;
  for m = 1:nm_test
    DATA = Con_ODE(f, L, u, [0; wait], X_TEST(:, m));
    nt = size(DATA, 1);
    C = DATA(nt, nx + 1);   % cost
    C_sum = C_sum + C;
  end  % m
  C_avg(1) = UCF*C_sum/nm_test
  figure(1)
  subplot(3, 1, 1), plot(UCF*(1:nt), DATA(:, 1:2))
  subplot(3, 1, 2), plot(UCF*(1:nt), DATA(:, 3), 'm')
  subplot(3, 1, 3), plot(UCF*(1:nt), DATA(:, 4))

  for r = 1:nr

  % Learn Jx*G
  YY = zeros(nf); YR = zeros(1, nf);
  for m = 1:np
    x = 2*rand(nx, 1) - 1; x = b*sin(0.5*pi*x);
    DATA = Jx_ODE(F, L_, DF_Dx, DL_Dx, [0 wait], x);
    nt = size(DATA, 1);
    for t = 1:nt
      x = DATA(t, 1:2)'; y = dz_dx(x)*G;
      Jx = DATA(t, 3:4); ref = Jx*G;
      YY = YY + y*y'; YR = YR + ref*y';
    end
  end  % m
  w = YR*pinv(YY);

  % Compute costs of new controller
  w_prev = w;
  u = @(x) -tanh(0.5*G'*dz_dx(x)'*w_prev');
  Du_Dx = @Dun_Dx_AK52_246;
  F = @(x) f(x, u(x));
  DF_Dx = @(x) df_dx(x) + df_du*Du_Dx(x);
  L_ = @(x) L(x, u(x));
  DL_Dx = @(x) dL_dx(x) + dL_du(u(x))*Du_Dx(x);
  C_sum = 0;
  for m = 1:nm_test
    DATA = Con_ODE(f, L, u, [0; wait], X_TEST(:, m));
    nt = size(DATA, 1); C = DATA(nt, nx + 1);  % cost
    C_sum = C_sum + C;
  end  % m
  C_avg(r + 1) = UCF*C_sum/nm_test
  figure(r + 1)
  subplot(3, 1, 1), plot(UCF*(1:nt), DATA(:, 1:2))
  subplot(3, 1, 2), plot(UCF*(1:nt), DATA(:, 3), 'm')
  subplot(3, 1, 3), plot(UCF*(1:nt), DATA(:, 4))

  end  % r
```

```matlab
function DATA = Con_ODE(f, L, u, T, xi)
% Simple numerical integrator.
% T in ms.

x = xi; u1 = u(x); nx = size(x, 1); nu = size(u1, 1);
C = 0; t_span = T(2) - T(1); DATA = zeros(nx + 1 + nu);
for t = 1:t_span
  u1 = u(x); L1 = L(x, u1); Dx1 = f(x, u1); x2 = x + Dx1;
  DATA(t, :) = [x', C, u1'];   % each row holds data for 1 time step
  u2 = u(x2); L2 = L(x2, u2); Dx2 = f(x2, u2); x3 = x + (Dx1 + Dx2)/4;
  u3 = u(x3); L3 = L(x3, u3); Dx3 = f(x3, u3); x = x + (Dx1 + 4*Dx3 + Dx2)/6;
  C = C + (L1 + 4*L3 + L2)/6;  % ignore UCF here, then correct in main file
end  % t

function Du_Dx = Du1_Dx_AK52(x)
% Derivative of initial controller wrt state x.
u = -5*x(1) - 3*x(2);
if abs(u) > 1, Du_Dx = zeros(1, 2); else Du_Dx = [-5, -3]; end
end

function Dun_Dx = Dun_Dx_AK52_246(x)
global G w_prev
d2z_dx2_G = G(2)*[0, 0;
    1, 0;
    0, 2;
    0, 0;
    3*x(1)^2, 0;
    2*x(1)*2*x(2), x(1)^2*2;
    3*x(2)^2, x(1)*6*x(2)
    0, 12*x(2)^2;
    0, 0;
    5*x(1)^4, 0;
    4*x(1)^3*2*x(2), x(1)^4*2;
    3*x(1)^2*3*x(2)^2, x(1)^3*6*x(2);
    2*x(1)*4*x(2)^3, x(1)^2*12*x(2)^2;
    5*x(2)^4, x(1)*20*x(2)^3;
    0, 30*x(2)^4];
Dun_Dx = -(1 - tanh(0.5*G'*dz_dx_AK52_246(x)'*w_prev')^2)*0.5*w_prev*d2z_dx2_G;
end

function dz_dx = dz_dx_AK52_246(x)
% Feature gradient: derivatives of 2nd-, 4th-, 6th- & 8th-deg
% monomials of 2 scalars.
dz_dx = [2*x(1), 0;
         x(2), x(1);
         0, 2*x(2);
         4*x(1)^3, 0;
```

```
            3*x(1)^2*x(2), x(1)^3;
            2*x(1)*x(2)^2, x(1)^2*2*x(2);
            x(2)^3, x(1)*3*x(2)^2;
            0, 4*x(2)^3;
            6*x(1)^5, 0;
            5*x(1)^4*x(2), x(1)^5;
            4*x(1)^3*x(2)^2, x(1)^4*2*x(2);
            3*x(1)^2*x(2)^3, x(1)^3*3*x(2)^2;
            2*x(1)*x(2)^4, x(1)^2*4*x(2)^3;
            x(2)^5, x(1)*5*x(2)^4;
            0, 6*x(2)^5];
end

function DATA = Jx_ODE(F, L_, DF_Dx, DL_Dx, T, xi)
% Forward- and backward-in-time integrations.

L_fin = 1e-6;

% Perform movement
x = xi; nx = size(x, 1); nt = T(2) - T(1);
X = zeros(nt, nx);
for t = 1:nt - 1
  X(t, :) = x';
  Dx1 = F(x); x2 = x + Dx1;
  Dx2 = F(x2); x3 = x + (Dx1 + Dx2)/4;
  Dx3 = F(x3); x = x + (Dx1 + 4*Dx3 + Dx2)/6;
  L1 = L_(x); if L1 < L_fin, break; end
end   % t
nt = t + 1;
X(nt, :) = x';
if L1 > L_fin, L1, end

% Calculate backward in time to find Jx
DATA = zeros(nt, 2*nx);
Jx = zeros(1, nx);
for t = nt:-1:2
  DATA(t, :) = [x' Jx];
  DJx1 = DL_Dx(x) + Jx*DF_Dx(x); Jx2 = Jx + DJx1;
  Dx1 = -F(x); x2 = x + Dx1;
  DJx2 = DL_Dx(x2) + Jx2*DF_Dx(x2); Jx3 = Jx + (DJx1 + DJx2)/4;
  Dx2 = -F(x2); x3 = x + (Dx1 + Dx2)/4;
  DJx3 = DL_Dx(x3) + Jx3*DF_Dx(x3); Jx = Jx + (DJx1 + 4*DJx3 + DJx2)/6;
  Dx3 = -F(x3); x = x + (Dx1 + 4*Dx3 + Dx2)/6;
  x = X(t - 1, :)';  % breadcrumb
  Dx = F(x); Dxip = Dx'*Dx; L1 = L_(x); Jx = Jx - ((Jx*Dx + L1)/Dxip)*Dx';   % GHJB
end
DATA(1, :) = [x' Jx];
```

```
function u = u1_AK52(x)
% Initial controller.
u = -5*x(1) - 3*x(2);
ao = 1 - 1000*eps;
if abs(u) > 1, u = ao*sign(u); end
end

function z = z_AK52_246(x)
% Features: 2nd-, 4th- & 6th-deg monomials of 2 scalars.
z = [x(1)^2; x(1)*x(2); x(2)^2;
     x(1)^4; x(1)^3*x(2); x(1)^2*x(2)^2; x(1)*x(2)^3; x(2)^4;
     x(1)^6; x(1)^5*x(2); x(1)^4*x(2)^2; x(1)^3*x(2)^3;
     x(1)^2*x(2)^4; x(1)*x(2)^5; x(2)^6];
end
```